\documentclass[12pt]{amsart}
\usepackage{amsmath,amsthm,amssymb}
\author{Boaz Tsaban}
\thanks{This paper is based on my M.Sc.\ thesis work supervised by Martin Goldstern.}
\address{Department of Mathematics and Computer Science, Bar-Ilan University,
52900 Ramat-Gan, Israel.}
\email{tsaban@macs.biu.ac.il}
\urladdr{http://www.cs.biu.ac.il/\~{}tsaban}

\title{A topological interpretation of the tower number}
\dedicatory{This paper is dedicated to Prof.\ Rabbi Haim Judah}

\newcommand{\arx}[1]{\texttt{http://arxiv.org/abs/#1}}
\newcommand{\ft}{\mathfrak{t}}
\long\def\forget#1\forgotten{}

\newcommand{\ot}{\stackrel}
\newcommand{\Cal}{\mathcal}

\newcommand{\CA}{{\Cal A}}

\newcommand{\F}{{\Cal F}}
\newcommand{\G}{{\Cal G}}
\newcommand{\J}{{\Cal J}}

\newcommand{\MEN}{M\!E\!N}

\newcommand{\R}{\mathbb R}

\newcommand{\Tau}{{\Cal T}}
\newcommand{\U}{\bigcup}

\renewcommand{\b}{\mathfrak{b}}
\renewcommand{\c}{\mathfrak{c}}
\renewcommand{\d}{\mathfrak{d}}

\newcommand{\g}{\gamma}
\renewcommand{\i}{\item}
\renewcommand{\k}{\kappa}
\newcommand{\oo}{\infty}
\newcommand{\p}{\mathfrak{p}}
\newcommand{\s}{\sigma}
\newcommand{\uw}{{\uparrow\omega}}
\newcommand{\w}{\omega}

\newcommand{\nin}{\not\in}
\newcommand{\cat}{\hat{\ }}
\newcommand{\sbst}{\subseteq}
\newcommand{\spst}{\supseteq}
\newcommand{\sm}{\setminus}
\newcommand{\as}{\subseteq^*}
\renewcommand{\|}{\upharpoonright}

\renewcommand{\)}{\right)}
\newcommand{\<}{\langle}
\renewcommand{\>}{\rangle}
\newcommand{\A}{\forall}

\newcommand{\op}{\operatorname}

\newcommand{\add}{\op{add}}

\newcommand{\non}{\op{non}}

\long\def\note#1\endnote%
{\ifNote{\par\medskip\tt Note: #1}\medskip\par
                         \else{}
                         \fi}

\newcommand{\impl}{\to}
\renewcommand{\t}{\tilde}

\newtheorem{theorem}{Theorem}[section]
\newtheorem{cor2}[theorem]{Corollary}

\newtheorem{lem2}[theorem]{Lemma}

\newtheorem{prob2}[theorem]{Problem}

\theoremstyle{remark}

\theoremstyle{definition}
\newtheorem{rem2}[theorem]{Remark}
\newtheorem{Qn}[theorem]{Question}

\newcommand{\be}{\begin{enumerate}}
\newcommand{\ee}{\end{enumerate}}
\newcommand{\bpf}{\begin{proof}[\sc Proof. ]}
\newcommand{\epf}{\end{proof}}

\def\Clubsuit{{\footnotesize ($\scriptstyle{\clubsuit}$)}}

\begin{document}
\maketitle

\begin{abstract}
Hurewicz
 found connections between some topological notions and the
combinatorial cardinals $\b$ and $\d$.
Rec\l{}aw
gave topological meaning to the definition of the cardinal $\p$.
We extend the picture
with a topological interpretation of the cardinal $\ft$.
We compare our notion to the one related to $\p$, and to some other
classical notions.
This sheds new light on the famous open problem whether $\p=\ft$.
\end{abstract}

\section{Introduction}

Cardinals associated with infinitary combinatorics play an important
role in set theory. Some earlier works
(\cite{HURE2}, \cite{RECLAW}, \cite{JUBAR}, \cite{PAWREC}, and
\cite{OPEN}) have pointed out a strong connection
between these cardinals and classes of spaces having certain
topological
properties. In this paper, we continue this line of research in a
way which enables us to give a topological
meaning to an open problem from infinitary combinatorics.

\subsection{Preliminaries}
Let $\w=\{0,1,2,\dots\}$ and $2=\{0,1\}$ be the usual discrete
spaces. $\w^\w$ and $2^\w$ are equipped with
the product topology.
Identify $2^\w$ with $P(\w)$ via characteristic functions.
$[\w]^\w$ is the set of infinite elements of $P(\w)$, with
$O_n=\{a:n\in a\}$ and $O_{\lnot n}=\{a:n\nin a\}=O_n^c$
($n\in\w$) as a
clopen subbase.

For $a, b\sbst\w$, $a\as b$ if $a\sm b$ is finite.
 $X\sbst [\w]^\w$ is \emph{centered} if every finite
$F\sbst X$ has an infinite intersection. $a\in [\w]^\w$ is an
\emph{pseudo-intersection} of $X$ if it is infinite, and
for all $b\in X$, $a\as b$.
$X\sbst[\w]^\w$ is a \emph{power} if it is centered, but has no pseudo-intersection.
$\p$ is the minimal size of a power.
$X\sbst [\w]^\w$ is a \emph{tower} if it is linearly ordered by
$\as$, and
has no pseudo-intersection{}.%
\footnote{
Note that, unlike the customary definition, we do not demand that a
tower is well-ordered.
However, by \cite[Theorem 3.7]{vD}, this does not change the value
of $\ft$.
}
$\ft$ is the minimal size of a tower.

$\le^*$ is the partial order defined on $\w^\w$ by eventual
dominance
($f\le^* g$ iff $\A^\oo n\ \(f(n)\le g(n)\)$~).
$\b$ is the minimal size of an unbounded family, and $\d$ is the
minimal size of a dominating family, with
respect to $\le^*$.

\subsubsection*{The Main Problem}
Let $\c$ denote the size of the continuum. The following holds.

{
\begin{theorem}[{\cite[Theorem 3.1.a]{vD}}]
 $\p\le\ft\le \b\le\d\le\c$.
\end{theorem}
}

For each pair except $\p$ and $\ft$, it is well known that a strict
inequality is consistent.

\begin{prob2}\label{p-t}
Is $\p < \ft$ consistent with $ZFC$ ?
\end{prob2}

This is one of the major and oldest problems of infinitary
combinatorics. Allusions for this problem can be
found in Rothberger's works (see, e.g., \cite[Lemma 7]{ROTH2}).
It is only known that $\p=\w_1\impl \ft=\w_1$
\cite[Theorem 3.1.b]{vD}, hence also $\ft\le\w_2\impl\p=\ft$.

\forget
Hurewicz \cite{HURE2} and Rec\l{}aw \cite{RECLAW} have compared some
topological notions to combinatorial notions
associated with the cardinals $\p$, $\b$ and $\d$.
The following section will describe the topological notion associated
with $\p$.
This will serve as a motivation for our presentation of a new
topological notion
in the third section, which will turn out to be associated to the
definition of $\ft$, and hence
yield a topological characterization of $\ft$.
We will also compare our notion to the one used by Rec\l{}aw.
In the fourth section, we will compare this notion to other classical
 notions.
The last section will introduce and study a closely related notion.
We hope that a study of our characterization
might lead to a better understanding of the above open problem, and
--- perhaps --- to its solution.
\forgotten

\subsection{$\g$-spaces}
Throughout this paper, by \emph{space} we mean a zero-dimensional,
separable, met\-riz\-able topological space.
As any such space is homeomorphic to a subset of the irrationals,
our results can be thought
of as dealing with sets of reals.

The definition of a $\g$-space is due to Gerlits and Nagy \cite{GN}.
Let $X$ be a space.
A collection of open sets $\G$ is an \emph{$\w$-cover} of $X$ if for
every finite
$F\sbst X$ there is a $G\in\G$ such that $F\sbst G$.
$\<G_n:n<\w\>$ is a \emph{$\g$-sequence} for $X$ if
$\A x\in X\A^\oo n\(x\in G_n\)$.
An open cover $\G$ of $X$ is a \emph{$\g$-cover} of $X$ if it
contains a $\g$-sequence for $X$.%
\footnote{In more recent papers, the term $\g$-cover is used for
what we call $\g$-sequence in this paper. BT}
Clearly every $\g$-cover is an $\w$-cover.
$X$ is a \emph{$\g$-space} if every $\w$-cover of $X$ is a
$\g$-cover of $X$.
For convenience, we may assume that the $\w$-covers $\G$ of $X$
are countable and clopen
(replacing each element from $\G$ by all finite unions of basic
clopen sets contained in it), and
that for all finite $F\sbst X$, there are \emph{infinitely many}
$G\in\G$ with $F\sbst G$.
(Using a partition of $\w$ into $\w$ many infinite sets, one can
create a sequence consisting
of $\w$ copies of the original cover $\G$. The resulted sequence
has the required property.)
And so we will, from now on.

Let $\Gamma$ denote the collection of all $\g$-spaces.
Rec\l{}aw has given an elegant characterization of $\g$-spaces.

{
\begin{theorem}[Rec\l{}aw {\cite[Theorem 3.2]{RECLAW}}]%
\label{gamma-to-P}
 $X$ is a $\g$-space iff no continuous image%
\footnote{
A continuous image of $X$ is the image of a continuous function
with domain $X$.
}
of $X$ in $[\w]^\w$ is a power.
\end{theorem}
}

This gives an alternative proof to the following. For a family $\F$
of spaces, let
$$\non(\F)=\min\{|X| : X\mathrm{\ is\ a\ space,\ } X\nin \F\}$$

\begin{cor2}[Galvin, Miller, Taylor {\cite[p.~146]{GM}}]%
\label{gamma-p}
$\non(\Gamma)=\p$.
\end{cor2}

\section{The Tower of $\tau$}

Let $X$ be a topological space.
An open cover $\G$ of $X$ is \emph{$T_1$} if for all distinct
$x,y\in X$ there is a $G\in\G$ such that
$x\in G, y\nin G$. Therefore,
 $\G$ is \emph{not $T_1$} iff there are distinct $x,y\in X$ such
that for all $G\in\G$,
 $x\in G\to y\in G$. Strengthening the ``not $T_1$'' property
demanding that
the above holds for \emph{all} $x,y\in X$ would trivialise $\G$ to
be $\{X\}$, or
$\{X,\emptyset\}$.
We therefore compensate by means of a ``modulo finite'' restriction.
\forget
loosen the restriction to be, that for all $x,y$, either
 $x\in G\to y\in G$ except for finitely many $G\in\G$, or
 $y\in G\to x\in G$ except for finitely many $G\in\G$. We put
this
definition more precisely.
\forgotten

\newcommand{\ldsto}{\mathrel{\ot\G\rightsquigarrow}}
For $\G=\<G_n:n<\w\>$, we write $x\ldsto y$ for
$$\A^\oo n\(x\in G_n\to y\in G_n\).$$
$\G$ is a \emph{$\tau$-sequence} for $X$ if
\be
\item $\G$ is a \emph{large cover}; i.e. every element of $X$ is
covered by infinitely many
elements of $\G$,%
\footnote{
This requirement was added in order to avoid trivial cases.
}
and
\item for all $x,y\in X$, either $x\ldsto y$, or $y\ldsto x$.
\ee

An open cover $\J$ of $X$ is a \emph{$\tau$-cover} of $X$ if it
contains a $\tau$-sequence for $X$.%
\footnote{In more recent papers, the term $\tau$-cover is used for
what we call $\tau$-sequence in this paper. BT}
It is easy to see that every $\g$-cover is a $\tau$-cover, and
every $\tau$-cover is an $\w$-cover.

$X$ is a \emph{$\tau$-space} if every \emph{clopen} $\tau$-cover of $X$
is a $\g$-cover of $X$.
Equivalently, if every countable clopen $\tau$-cover of $X$ is a
$\g$-cover of $X$.
Let $\Cal T$ denote the collection of all $\tau$-spaces.

\begin{cor2}\label{Gamma-T}
$\Gamma\sbst \Cal T$.
\end{cor2}

We wish to transfer our covering notions into $[\w]^\w$, in order
to obtain their
combinatorial versions.
In Rec\l{}aw's proof of Theorem \ref{gamma-to-P},
the following natural function $h=h_\G$ is considered.
Given a countable large
cover $\G=\{G_n:n\in\w\}$ of $X$, define $h:X\to [\w]^\w$
by
$h(x)=\{n:x\in G_n\}.$
Now, let us see what $h$ does to our topological notions.
Assume that $\G=\<G_n:n<\w\>$ is an $\w$-cover of $X$. Then
for all finite $F\sbst X$, $F$ is a subset of infinitely many
$G_n$'s.
That is, $n\in \cap h[F]$ for infinitely many $n$'s. This means that
$h[X]$ is centered.

Next, assume that $\G$ is a $\g$-sequence for $X$. Then
$\A x\A^\oo n\ (x\in G_n)$.
That is, $\A x\A^\oo n\ (n\in h(x))$, or
$\w$ is a pseudo-intersection of $h[X]$. Therefore, $\G$ is a $\g$-\emph{cover}
of $X$
iff the associated $h[X]$ has a pseudo-intersection.

Finally, a large cover $\G$ is a $\tau$-sequence for $X$ iff
for all $x,y\in X$, either $x\ldsto y$, or $y\ldsto x$.
Now, $a\ldsto b$ iff $\A^\oo n\ (n \in h(a)\to n\in h(b))$ iff
$h(a)\as h(b)$.
Therefore, $h[X]$ is linearly ordered by $\as$.

We have showed the following.
\begin{lem2}\label{notions}
Assume that $\G$ is a countable large cover of $X$.
\be
\item $\G$ is an $\w$-cover of $X$ iff $h_\G[X]$ is centered.
\item $\G$ is a $\g$-cover of $X$ iff $h_\G[X]$ has an
almost-intersection.
\item $\G$ is a $\tau$-sequence for $X$ iff $h_\G[X]$ is
linearly ordered by $\as$.
\ee
\end{lem2}

Note that if $\G$ is a \emph{clopen} cover, then $h=h_\G$
is continuous, since for all $n$, $h^{-1}[O_n]=G_n$, and
$h^{-1}[O_{\lnot n}]=G_n^c$.
Therefore,
\ref{notions}(1) and \ref{notions}(2) yield Rec\l{}aw's Theorem
\ref{gamma-to-P}, and
\ref{notions}(2) and \ref{notions}(3) yield the following.

{
\begin{theorem}\label{tauimage}
$X$ is a $\tau$-space iff no continuous image of $X$ in
$[\w]^\w$ is a tower.
\end{theorem}
}
\begin{proof}[\sc Proof. ]
($\Leftarrow$) Assume that $\J$ is a clopen $\tau$-cover of $X$
and let $\G\sbst\J$ be a $\tau$-sequence for $X$.
Then by Lemma \ref{notions}(3), $h_\G[X]$ is linearly ordered by
$\as$.
As $h_\G$ is continuous, $h_\G[X]$ cannot be a tower, and hence
has a pseudo-intersection.
Applying Lemma \ref{notions}(2), we get that $\G$ is a $\g$-cover
of $X$, and hence so is $\J$.

($\Rightarrow$) Assume that $X$ is a $\tau$-space, and
$f:X\to[\w]^\w$ is continuous.
Assume that $f[X]$ is linearly ordered by $\as$. Then
$\<O_n:n<\w\>$ is a clopen $\tau$-sequence for $f[X]$.
Therefore $\G=\<f^{-1}[O_n] : n\in\w\>$ is a clopen
$\tau$-sequence for $X$; hence a $\g$-cover of $X$.
By \ref{notions}(2) $h_\G[X]$
has a pseudo-intersection; hence is not a tower.
But $h_\G=f$, as for all $x\in X$,
$$n\in h_\G(x)\iff x\in f^{-1}[O_n]\iff f(x)\in O_n \iff n\in f(x).$$
Therefore, $f[X]$ is not a tower.
\end{proof}

The reader might have noticed that in the last proof we have
indirectly used the following lemma.

\begin{lem2}\label{cont-T}
Every continuous image of a $\tau$-space is a $\tau$-space.
\end{lem2}

\begin{proof}[\sc Proof. ]
A continuous preimage of a clopen $\tau$-cover is a clopen
$\tau$-cover.
\end{proof}

We get a topological characterzation of $\ft$.

\begin{cor2}
\label{T-t}
$\non(\Tau)=\ft$.
\end{cor2}

For a family $\F$ of spaces, let
$$\add(\F)=\min\{|\CA| : \CA\sbst\F\ \&\ \U\CA\nin\F\}.$$

\begin{theorem}\label{addtau}
$\add(\Tau) = \ft$.
\end{theorem}
\begin{proof}[\sc Proof. ]
We need the following lemma.
\begin{lem2}\label{pseudo-intersectionpieces}
Assume that $X\sbst [\w]^\w$ is linearly ordered by $\as$, and
for some $\k<\ft$, $X=\U_{i<\k}X_i$ where each $X_i$ has a pseudo-intersection.
Then $X$ has a pseudo-intersection.
\end{lem2}
\begin{proof}[\sc Proof. ]
If for each $i<\k$ $X_i$ has a pseudo-intersection $x_i\in X$, then a pseudo-intersection of
$\{x_i : i<\k\}$
is also a pseudo-intersection of $X$.
Otherwise, there exists $i<\k$ such that $X_i$ has no pseudo-intersection
$x\in X$.
That is, for all $x\in X$ there exists a $y\in X_i$ such that
$x\not\as y$; thus
$y\as x$.
Therefore, a pseudo-intersection of $X_i$ is also a pseudo-intersection of $X$.
\end{proof}
Now we can use Theorem \ref{tauimage}.
Assume that $X=\U_{i<\k<\ft}X_i$, where each $X_i$ is a
$\tau$-space.
If $f:X\to [\w]^\w$ is continuous and $f[X]$ is linearly ordered
by $\as$,
then each $f[X_i]$ has a pseudo-intersection and therefore by the lemma, $f[X]$
has a pseudo-intersection.
Therefore, $X$ is a $\tau$-space.
\end{proof}

A similar result cannot be obtained for $\g$-spaces. $CH$ implies
that $\g$-spaces
are not even closed under taking \emph{finite} unoions. We will use an
argument of Galvin and Miller
\cite[p.~151]{GM} to show this.

\begin{theorem}[Brendle {\cite[Theorem 4.1]{JORG}}]%
\label{Brendlespace}
$CH\to$ there is a $\g$-space of size $\c(=\w_1)$ all of whose
subsets are $\g$-spaces.
\end{theorem}

For $X\sbst [0, 1]$, let $X+1 = \{x+1 : x\in X\}$.
\begin{theorem}[Galvin, Miller {\cite[Theorem 5]{GM}}]
Assume that $A\sbst X\sbst [0,1]$, and $(X\sm A)\cup (A+1)$
is a $\g$-space.
Then $A$ is $G_\delta$ and $F_\sigma$ in $X$.
\end{theorem}

Now, consider a subspace $A$ of Brendle's space $X$, such that $A$
is neither $G_\delta$
nor $F_\sigma$. Then the union of $X\sm A$ and $A+1$ (which
are both $\g$-spaces)
is not a $\g$-space.

\medskip

The $\g$-property is very strict. Gerlits and Nagy
\cite[Corollary 6]{GN} proved that
$\g$-spaces are $C''$.
In particular, it is consistent that all $\g$-spaces are countable.
However, large $\tau$-spaces do exist in $ZFC$.
In fact, we have the following.

\begin{theorem}[Shelah]\label{Shelah}
$2^\w$ is a $\tau$-space.
\end{theorem}

\begin{proof}[\sc Proof. ]
Towards a contradiction, assume that $f:2^\w\to [\w]^\w$ is
continuous such that $f[2^\w]$ is a tower.
Let
$$T =  \{ s \in 2^{<\w}:  f[\ [s]\ ] \mathrm{\ has\ no\ pseudo-intersection}\}.$$
$T$ is a perfect tree.
Assume that for some $s\in T$,
there are no incomparable extensions $s_0,s_1$ such that both
$f[~[s_0]~]$ and $f[~[s_1]~]$ have no
pseudo-intersection. Then for all $\t s$ extending $s$, at least one of
$f[~[\t s\cat\<0\>]~]$ and $f[~[\t s\cat\<1\>]~]$ has
a pseudo-intersection.
Let $\s\in 2^\w$ extend $s$ such that for all $n\ge |s|$,
$f[~[(\s \| n)\cat \<1-\s(n)\>]~]$ has a pseudo-intersection.
$f[~[s]~] = \U_{n<\w} f[~[(\s \| n)\cat \<1-\s(n)\>]~]
\cup \{f(\s)\}$ is a union of
$\w$ many sets having a pseudo-intersection, contradicting Lemma \ref{pseudo-intersectionpieces}.

We now show that $f[2^\w]$ cannot be linearly ordered by $\as$.
Define two branches $\beta$ and $\xi$ in $T$ as follows.   Start
with
incomparable $b_0,c_0\in T$.  Pick $x_0 \in [b_0]$. As
$f[~[c_0]~]$
has no pseudo-intersection, we can find a $y_0\in [c_0]$ such that
$f(x_0)\not\as f(y_0)$.
Choose an $n_0\in f(x_0)\sm f(y_0)$.
Since $f$ is continuous, we can find $b_1$, an initial segment of
$x_0$,
such that $f[~[b_1]~]\sbst O_{n_0}$.
Similarly, find $c_1$, an initial segment of $y_0$, such that
$f[~[c_1]~]\sbst O_{\lnot n_0}$.

Now we reverse the roles, and find $x_1\in [b_1]$,
$y_1\in [c_1]$, $n_1> n_0$
such that $n_1\in f(y_1)\sm f(x_1)$.
Then we take $b_2$ and $c_2$, initial segments of $y_1$ and
$x_1$ respectively, such that
$f[~[b_2]~]\sbst O_{\lnot n_1}$ and $f[~[c_2]~]\sbst O_{n_1}$.

We continue by induction.
Finally, let $\beta=\U_i b_i=\lim_i x_i$, and
$\xi = \U_i c_i=\lim_i y_i$.
Since $f$ is continuous, the sets $\{ n_{2k} : k\in\w \}$ and
$\{n_{2k+1} : k\in\w \}$ witness that
neither $f(\beta) \as f(\xi)$ nor $f(\xi) \as f(\beta)$.
\end{proof}

This theorem implies that the inclusion in Corollary \ref{Gamma-T}
is proper. We will
modify it to get a large class of $\tau$-spaces which are not
$\g$-spaces.

\begin{theorem}\label{w^w}
$\w^\w$ is a $\tau$-space.
\end{theorem}

\begin{proof}[\sc Proof. ]
Identify $\w^\w$ with $2^\w\sm F$ (where $F$ are the
eventually zero sequences), and
work in $2^\w\sm F$ instead of $2^\w$.
\be
\i In the proof that $T$ is perfect, we need not care whether
$\s\in 2^\w\sm F$ or not.
\i When choosing the initial segment $b_{i+1}$ of $x_i$,
use the fact that $x_i\nin F$ to make sure that $b_{i+1}$ ends
with a ``$1$'' (a similar treatment for $c_{i+1}$).
This will make $\beta$ and $\xi$ belong to $2^\w\sm F$.
\ee
\end{proof}

\begin{cor2}\label{analytic}
Every analytic set of reals is a $\tau$-space.
\end{cor2}
\begin{proof}[\sc Proof. ]
Every analytic set of reals is a continuous image of $\w^\w$.
\end{proof}

\begin{rem2}\label{suslin} \
\be
\i One cannot prove in $ZFC$ that all projective sets of reals are
$\tau$-spaces.
Since the reals have a projective well-ordering in the constructible
 universe $L$, a straightforward
inductive construction will yield a projective tower.

\i Due to a theorem of Suslin (see, e.g.,
\cite[Corollary 2C.3]{MOSCHO}),
every uncountable analytic set contains a perfect set, and hence
is not a $\g$-space. (It is not even strongly null.)
\ee
\end{rem2}

As in the case of $\g$-spaces \cite[p.~147]{GM},
the property of being a $\tau$-space need not be hereditary for
subspaces of the same size.

We will work in $P(\w)$.

\begin{theorem}\label{tau-not-hereditary}
$\ft=\c\to$ there is a space $X\sbst [\w]^\w$  s.t.
\be
\i $|X|=\c$,
\i $X\cup [\w]^{<\w}$ is a $\tau$-space, and
\i $X$ is not a $\tau$-space.
\ee
\end{theorem}
\begin{proof}[\sc Proof. ]
First, note that (1) follows from (2) and (3), using Corollary
\ref{T-t}.

We will use a modification of the Galvin-Miller construction (see
\cite[Theorem 1]{GM}).
For $y\in [\w]^\w$, define $y^*=\{x : x\as y\}$.
We need the following lemma.

\begin{lem2}[{Galvin, Miller \cite[Lemma 1.2]{GM}}]
Assume that $\G$ is an open $\w$-cover of $[\w]^{<\w}$.
Then for all $x\in [\w]^\w$ there exists a $y\in [x]^\w$ such
that $\G$ $\g$-covers $y^*$.
\end{lem2}

Let $\<\G_i:i<\c\>$ enumerate all countable families of clopen
sets in $P(\w)$, and let
$\<y_i:i<\c\>$ enumerate all elements $y\in [\w]^\w$ such that
both $y$ and $\w\sm y$ are infinite.

Construct, by induction, $\<x_i:i<\c\>\sbst [\w]^\w$ such that
$i<j\to x_j\as x_i$.
For a limit $i$, use $i<\ft$ to get $x_i$. For successor $i=k+1$,
$x_i$ is constructed as follows.
\begin{description}
\item[Case 1] $\G_k$ is a $\tau$-cover of
$B_k=\{x_j:j\le k\}\cup [\w]^{<\w}$.
By Theorem \ref{T-t}, as $|B_k|<\ft$, $\G_k$ is a $\g$-cover of
$B_k$.
In particular, $\G_k$ $\g$-covers $[\w]^{<\w}$.
By the lemma, there exists an $x_{k+1}\in [x_k]^\w$ such that
$\G_k$ $\g$-covers $x_{k+1}^*$.
\item[Case 2] $\G_k$ is not a $\tau$-cover of
$\{x_j:j<k\}\cup [\w]^{<\w}$. Since this case
is not interesting, we may take $x_{k+1}=x_k$.
\end{description}

After $x_i$ is chosen (either for limit or successor $i$), modify it
as follows. If $x_i\as y_i$, leave it as is.
Otherwise, replace it by $x_i\sm y_i$.
This does the construction.

Define $X=\{x_i:i<\c\}$. Then $X\cup [\w]^{<\w}$ is a
$\tau$-space.
By the construction, if $\G_k$ is a $\tau$-cover of
$X\cup [\w]^{<\w}$, then it $\g$-covers
$\{x_j:j\le k\}\cup x_{k+1}^*$. But
$$X\cup [\w]^{<\w}\sbst\{x_j:j\le k\}\cup x_{k+1}^*.$$
This does (2).

(3) $X$ is a tower. Let $a\in [\w]^\w$. We will show that $a$ is
not a pseudo-intersection of $X$.
Take $a_0\sbst a$ such that both $a_0$ and $\w\sm a_0$ are
infinite. Now, some $x_i$ satisfies either
$x_i\as a_0$, or $x_i\as \w\sm a_0$. Therefore,
$a\not\as x_i$.
By Lemma \ref{cont-T} (considering the identity function on
$[\w]^\w$), $X$ is not a $\tau$-space.
\end{proof}

\begin{cor2}\label{not-Borel}
$\ft=\c\to$ $\tau$-spaces are not closed under Borel images.
\end{cor2}
\begin{proof}[\sc Proof. ]
Let $X$ be given by the theorem.
Consider any function $f: X\cup [\w]^{<\w}\to [\w]^\w$ such
that $f\| X$ is the identity function, and
$f[[\w]^{<\w}]\sbst X$.
As $[\w]^{<\w}$ is countable, $f$ is Borel. $X\cup [\w]^{<\w}$
is a $\tau$-space, but $X$, its Borel image,
is not a $\tau$-space.
\end{proof}

\section{Comparing $\tau$-spaces to Other Classical Classes}

\subsection*{Hurewicz and Menger}

We give Hurewicz' topological interpretations of $\b$ and $\d$.

$X$ has the \emph{Hurewicz property} if for every sequence of open
covers $\G_n$,
there is a sequence of finite $\t\G_n\sbst \G_n$ such that the
sets $\cup \t\G_n$
form a $\g$-cover of $X$.
$X$ has the \emph{Menger property} if for every sequence of open
covers $\G_n$,
there is a sequence of finite $\t\G_n\sbst \G_n$ such that the
sets $\cup \t\G_n$
cover $X$. Let $\Cal H$ and $\MEN$ denote the classes of spaces
having the Hurewicz and Menger properties,
respectively.
Clearly ${\Cal H}\sbst \MEN$.

\begin{theorem}[Hurewicz {\cite[\S 5]{HURE2}}]\label{b-d-chars}
Let $X$ be a space.
\be
\i $X$ has the Hurewicz property iff every continuous image of $X$
in $\w^\w$ is bounded.
In particular, $\non({\Cal H})=\b$.
\i $X$ has the Menger property iff every continuous image of $X$ in
$\w^\w$ is not dominating.
In particular, $\non(\MEN)=\d$.
\ee
\end{theorem}

We get that none of these two notions is provably comparable to
$\Tau$.

\begin{cor2} \
\be
\i $\Tau\not\sbst \MEN$, and
\i $\ft<\b\to{\Cal H}\not\sbst\Tau$
\ee
\end{cor2}
\begin{proof}[\sc Proof. ]
(1) By Theorem \ref{w^w}, $\w^\w\in\Tau$, and by Theorem
\ref{b-d-chars}(2),
$\w^\w\nin\MEN$.

(2) follows from Corollary \ref{T-t} and Theorem \ref{b-d-chars}(1).
\end{proof}

Indeed, $\tau$-spaces could be pretty
far from having the Menger property. According to a
theorem of Hurewicz \cite[Theorem 20]{HURE}, an analytic set of
reals having the Menger
property must be $F_\s$. Corollary \ref{analytic} could
be contrasted with this.
However, these classes need not be orthogonal.
Gerlits and Nagy \cite[p.~155]{GN} proved that, given a sequence of
$\w$-covers $\G_n$ of a $\g$-space $X$,
there exists a sequence $G_n\in\G_n$ such that
$\{G_n:n\in\w\}$ $\g$-covers $X$.
We therefore have the next assertion.
\begin{cor2}
$\Gamma\sbst {\Cal H}\cap \Tau$.
\end{cor2}

\subsection*{$\lambda$-spaces}

$X$ is a $\lambda$-\emph{space} if every countable subset of $X$ is
$G_{\delta}$.
Let $\Lambda$ denote the collection of $\lambda$-spaces.
$\lambda$-spaces are perfectly meager (see
\cite[Theorem 5.2]{MILLER}).
Therefore, by Remark \ref{suslin}(2), no uncountable analytic set
is a $\lambda$-space.
This again could be contrasted with Corollary \ref{analytic}.

On the other hand, we have the following.

\begin{theorem}
There is a $\lambda$-space of size $\ft$ which is not a
$\tau$-space.
\end{theorem}

Our theorem follows from the following two lemmas.

\begin{lem2}[{\cite[Theorem 1]{ROTH3}}]\label{lambda-b}
$\non(\Lambda)=\b$.
\end{lem2}

\begin{lem2}\label{b-lambda}
Every tower of size $\b$ is a $\lambda$-space.
\end{lem2}

\begin{proof}[\sc Proof. ]
We use the standard argument (see \cite[Theorem 9.1]{vD}). Before
getting started, note that
for all $y\in [\w]^\w$,
$y^*=\U_{s\in [\w]^{<\w}} \{x : x\sbst y\cup s\}$ is $F_\s$.

Assume that $X=\{x_i : i<\b\}$ is a tower with $i<j\to x_j\as x_i$.
For $\alpha <\b$, set $X_\alpha = \{x_i : i<\alpha \}$.
Then each $X_\alpha$ is $G_\delta$ in $X$. (Its complement in
$X$ is $F_\s$.)
Assume that $F\sbst X$ is countable. As $\b$ is regular,
there exists $\alpha<\b$ such that $F\sbst X_\alpha$.
As $|X_\alpha|<\b$, $X_\alpha$ is a $\lambda$-space.
Hence, $F$ is $G_\delta$ in $X_\alpha$; i.e., there is a
$G_\delta$ set $A\sbst X$ such that
$F=X_\alpha\cap A$. As $X_\alpha$ is also $G_\delta$, $F$
is $G_\delta$ in $X$.
\end{proof}

\forget
Therefore, $\ft <\b$ implies  that there is a $\lambda$-space of
size $\ft$ which is not a $\tau$-space
(consider a tower of size $\ft$).
The following theorem shows that this is already true in $ZFC$.
\forgotten

\forget
\begin{cor2}
If $\ft=\b$, or $\b < \d$,
then there is a a $\lambda$-space of size $\b$ which is not a
$\tau$-space.
\end{cor2}

\begin{proof}[\sc Proof. ]
In either case, consider a tower of size $\b$.
\end{proof}

\begin{rem2}
The increasing $\b$-sequence from the proof of Theorem
\ref{b<d-tower} was used by
Rothberger \cite[Theorem 4]{ROTH} as an example of a
$\lambda$-space which is not $\lambda'$. See also
\cite[Theorem 9.1]{vD}. We do not know whether there is a $ZFC$
example of a $\lambda$-space of size $\b$
which is not $\tau$.
\end{rem2}
\forgotten

With some set theoretic assumptions, we can have an example of size
$\b$. In fact,
our $\b$-example will have some additional properties related to
our study.
$X\sbst\R$ is a \emph{$\lambda'$-space} if for all countable
$F\sbst\R$, $X\cup F$ is a
$\lambda$-space. For $D\sbst\R$, $X$ is \emph{$\k$-concentrated}
on $D$ if for all open
$U\spst D$, $|X\sm U|<\k$.

Considering the proof that an $(\w_1, \w_1)$-gap is a
$\lambda'$-space (see \cite[p.~215]{MILLER}),
one might wonder whether our proof can be strengthened to make every
$\b$-tower $X$
a $\lambda'$-space.
In fact, following the steps of the proof carefully one gets that for
all countable $F\sbst [\w]^\w$,
$X\cup F$ is a $\lambda$-space. The problem is with
$[\w]^{<\w}$: If $X$, when viewed as a
subset of $\w^\w$, is unbounded, then $[\w]^{<\w}$ is not
$G_\delta$ in $X\cup [\w]^{<\w}$
\cite[Lemma 9.3]{vD}.

\begin{theorem}\label{specialX}
Assume that there exists a tower of size $\b$.
Then there is a space $X$ of size $\b$ such that
\be
\i $X$ is a $\lambda$-space,
\i $X$ is $\b$-concentrated on a countable set,
\i $X$ is not a $\lambda'$-space,
\i $X$ does not have the Hurewicz property, and
\i $X$ is not a $\tau$-space.
\ee
\end{theorem}

\begin{proof}[\sc Proof. ]
We work in $P(\w)$. Identify $[\w]^\w$ with $\w^\uw$, the
strictly increasing elements of $\w^\w$.
Let $X=\{x_i : i<\b\}\sbst \w^\uw$ be such that the following holds.

\smallskip
\centerline{\Clubsuit\quad It is unbounded,
$\le^*$-increasing, and has size $\b$.}
\smallskip

\noindent%
The existence of such a set follows from \cite[Theorem 3.3]{vD}.
Let $A=\{a_i : i<\b\}$ be a tower, and define $Y=\{y_i : i<\b\}$
as follows. For each $i<\b$,
let $h\in\w^\w$ bound $\{y_k : k<i\}\cup\{x_i\}$,
and take a $y_i\as a_i$ such that $h\le^*y_i$,
$y_i(n)=\min\{k\in a_i : y_i(n-1),h(n)<k \}.$
$Y$, like $X$, has the property \Clubsuit.
Rothberger \cite[Theorem 4]{ROTH} has proved that \Clubsuit{}
implies (1) and (3)
(see also \cite[Lemma 9.3]{vD}).
By an observation of Miller \cite[Theorem 5.7]{MILLER},
\Clubsuit{} implies that $Y$ is $\b$-concentrated on
$[\w]^{<\w}$.
By Theorem \ref{b-d-chars}(1), (4) is also satisfied.

(5) $Y$ is a tower: Any pseudo-intersection of $Y$ would also be a pseudo-intersection of $A$.
\end{proof}

Our theorem has a cute corollary.

\begin{cor2}
$\ft=\b\ \lor \b<\d\to$ there exists an $X$ as in Theorem
\ref{specialX}.
\end{cor2}

This follows from the following observation.

\begin{lem2}[{\cite[Theorem 1]{ZS}}]\label{b<d-tower}
$\b < \d\to$ there is a tower of size $\b$.
\end{lem2}

For completeness, we give a proof of this lemma.

\begin{proof}[\sc Proof. ]
Let $X\sbst \w^\uw$ have property \Clubsuit,
and let $h\in\w^\w$ witness that $X$ is not dominating.
For each $x\in X$ define $a_x\in [\w]^\w$ by
$a_x=\{n : x(n)<h(n)\}$. Then $\{a_x : x\in X\}$ is linearly
ordered by $\as$.
Assume that it has a pseudo-intersection $a$. Then $h'$ defined by
$h'(n)= h(\min\{k\in a : n\le k \})$
bounds $X$. A contradiction.
\end{proof}

Despite the large difference between $\tau$ and $\lambda$
spaces, these classes need not be orthogonal.
Their intersection
could contain a space of size $\c$:
By \cite[Theorem 2]{GM}, a $G_\delta$ $\g$-subspace of a space is
also an $F_\s$
subspace of that space.
Every co-countable subspace of Brendle's space (see Theorem
\ref{Brendlespace}) is $G_\delta$ and therefore $F_\s$.
Therefore, every countable subspace of Brendle's space is
$G_\delta$.

\begin{cor2}
$CH\to$ there is a  $\g$(in particular, $\tau$)-space of size
$\c$ which is also a $\lambda$-space.
\end{cor2}

\section{The Selection Principle $S_1$}

Unlike $\g$-spaces, $\tau$-spaces do not fit into the framework
defined in \cite{OPEN}.
We recall the basic definitions.

A space $X$ has property $S_1(x,y)$ ($x,y$ range over
$\{\w, \g, \tau,\dots\}$) if, given a \emph{sequence} of
$x$-covers $\G_n$, one can select from each $\G_n$ an element
$G_n$ such that $\{G_n:n\in\w\}$ is a $y$-cover.
As mentioned in section 3,
Gerlits and Nagy
proved that the $\g$-property is equivalent to the $S_1(\w,\g)$
property.
Using this notation, we have the following.
\begin{rem2}\label{inclusions}
$S_1(\w,\g)\sbst S_1(\tau,\g)\sbst S_1(\g,\g)$.
\end{rem2}
\begin{proof}[\sc Proof. ]
As noted in \S 3, every $\g$-cover is a $\tau$-cover, and every
$\tau$-cover is an $\w$-cover.
\end{proof}

Obviously, $S_1(\tau,\g)\sbst\Tau$.
By \cite[Theorem 2.3]{OPEN}, $2^\w$ does not belong to the class
$S_1(\g,\g)$,
and therefore not to the class $S_1(\tau,\g)$, either.

\begin{cor2}\label{no-tau-diagonal}
$S_1(\tau,\g)\neq\Tau$.
\end{cor2}

We now study the $S_1(\tau,\g)$ property.
Let us begin with saying that the $\tau$-covering notion fits
nicely into the framework of
\cite{OPEN}. (In fact, it suggests many interesting notions, but we
will stick to
$S_1(\tau,\g)$ in this paper.)
For example, it can be added to \cite[Theorem 3.1]{OPEN}.
In particular, we have the following.

\begin{theorem}
$S_1(\tau,\g)$ is closed under taking closed subsets and
continuous images.
\end{theorem}

There are more properties, which follow from Remark \ref{inclusions}
We quote some of them.

\begin{theorem}\label{properties} \
\be
\i {\rm (\cite[Corollary 5.6]{OPEN})}
Every element of $S_1(\tau,\g)$ is perfectly meager (i.e., has
meager intersection with every perfect set).
\i {\rm (\cite[Theorem 5.7]{OPEN})}
If $X\in S_1(\tau,\g)$, then for every $G_\delta$ set $G$
containing $X$, there exists an $F_\s$
set $F$ such that $X\sbst F\sbst G$.
\ee
\end{theorem}

\begin{rem2} If we omit the metrizability assumption on the
spaces, then $S_1(\tau,\g)$ is not closed under cartesian products,
nor under finite unions:
Todor\-{\v c}evi\'c \cite{TODOR} showed that there exist
nonmetrizable
$X,Y\in S_1(\w,\g)$ such that $X\cup Y\nin S_1(\g,\w)$. (In
fact,
he showed that they do not even have the Menger property.)
\end{rem2}

\begin{theorem}[Daniels {\cite[Lemma 9]{DANIELS}}]%
\label{gamma-powers}
$S_1(\w,\g)$ is closed under taking finite powers.
\end{theorem}

\begin{Qn}\label{powers}
Is $S_1(\tau,\g)$ closed under taking finite powers?
\end{Qn}

One can see that
if $\G$ is a $\tau$-cover of $X$, then $\{G^n:G\in\G\}$ is a
$\tau$-cover of $X^n$. But
this is not enough for answering this question.

\begin{theorem}
$\non(S_1(\tau,\g))=\ft$.
\end{theorem}
\begin{proof}[\sc Proof. ]
Assume that $|X|<\ft$ and let $\G_n$ be $\tau$-covers of $X$.
We wish to conclude that the $\G_n$'s are $\g$-covers of $X$.
Corollary \ref{T-t} is not enough for our purposes, since the
$\tau$-covers need not be clopen.
However, Theorem \ref{notions} gives the desired result.
Now, by \cite[Theorem 4.7]{OPEN}, $\non(S_1(\g,\g))=\b$.
As $|X|<\b$, $X\in S_1(\g,\g)$.
Therefore, one can extract a $\g$-cover of $X$ from the $\G_n$'s.
This proves $\ft\le\non(S_1(\tau,\g))$.

The other direction follows from the fact that
$S_1(\tau,\g)\sbst\Tau$, together with Corollary
\ref{T-t}.
\end{proof}

In particular, it is consistent that
$S_1(\tau,\g)\neq S_1(\g,\g)$.
\begin{cor2}
$\ft<\b\to S_1(\tau,\g)\neq S_1(\g,\g)$.
\end{cor2}

We therefore have the following.
\begin{Qn}\label{top-qn}
Does $S_1(\w,\g) = S_1(\tau,\g)$ ?
\end{Qn}

As the consistency of $\p<\ft$ would imply a negative answer, this
question seems to be closely related to the
main problem whether $\p=\ft$.%
\footnote{This problem is now solved, but tighter approximations to the
minimal tower problems are still open. See:
B.\ Tsaban, \emph{Selection principles and the minimal tower problem},
\arx{math.LO/0105045}}

\begin{rem2} \
\be
\i Due to Theorem \ref{properties}(1), the (in fact, any) Luzin set
used in \cite{OPEN} to distinguish
$S_1(\w,\w)$ from $S_1(\w,\g)$ will also distinguish it from
$S_1(\tau,\g)$.
\i By Theorem \ref{gamma-powers}, a negative answer to Question
\ref{powers} would imply a
negative answer to Question \ref{top-qn}.
\i Showing $S_1(\tau,\g)\not\sbst S_1(\w,\w)$ is consistent
would also yield a negative answer.
\ee
\end{rem2}

\section*{Acknowledgements}
{\small
I thank Martin Goldstern for supervising my thesis work (and the
writing of this paper).
I also thank Saharon Shelah for Theorem \ref{Shelah}, and
Arnold W.\ Miller, Ireneusz Rec\l{}aw, and the referees for making
some useful comments.
A special thanks is owed to Marion Scheepers, for very fruitful
discussions and suggestions which led to
\ref{tau-not-hereditary}, \ref{not-Borel}, and \S\S 3,4.
}

\end{document}